\documentclass[twoside,12pt]{article}
\usepackage{amssymb}
\usepackage{graphicx}
\usepackage{amsmath}
\usepackage{ams fonts}

\newtheorem{theorem}{Theorem}

\newtheorem{definition}{Definition}

\newtheorem{lemma}{Lemma}

\newtheorem{remark}{Remark}

\begin{document}
\begin{center}
\textbf{\Large Sasakian metric as a Ricci soliton\\
and related results}
\end{center}
Amalendu Ghosh$^1$ and Ramesh Sharma$^2$ 

\noindent
\textit{$^1$Department Of Mathematics, Krishnagar Government College, Krishnanagar 741101, 
West Bengal, India, E-mail: aghosh\_70@yahoo.com}\\ 
\noindent
\textit{$^2$Department Of Mathematics, University Of New Haven, West Haven, \\
CT 06516, USA, E-mail:rsharma@newhaven.edu}\\
\vskip.1in

\noindent
\textbf{Abstract: }We prove the following results: (i) A Sasakian metric as a non-trivial Ricci soliton is null $\eta$-Einstein, and expanding. Such a characterization permits to identify the Sasakian metric on the Heisenberg group $\mathcal{H}^{2n+1}$ as an explicit example of (non-trivial) Ricci soliton of such type. (ii) If an $\eta$-Einstein contact metric manifold $M$ has a vector field $V$ leaving the structure tensor and the scalar curvature invariant, then either $V$ is an infinitesimal automorphism, or $M$ is $D$-homothetically fixed $K$-contact.\\



\section{Introduction}
A Ricci soliton is a natural generalization of an Einstein metric, and is defined on a Riemannian manifold ($M,g$) by 
\begin{equation} \label{1}
(\pounds _{V}g)(X,Y) + 2Ric(X,Y) + 2\lambda g(X,Y) = 0
\end{equation}
where $\pounds_{V}g $ denotes the Lie derivative of $g$ along a vector field $V$, $\lambda $ a constant, and arbitrary vector fields $X,Y$ on $M$. The Ricci soliton is said to be shrinking, steady, and expanding accordingly as $\lambda$ is negative, zero, and positive respectively. Actually, a Ricci soliton is a generalized fixed point of Hamilton's Ricci flow \cite{Hamilton}: $\frac{\partial}{\partial t}g_{ij} = -2R_{ij}$, viewed as a dynamical system on the space of Riemannian metrics modulo diffeomorphisms and scalings. For details, see Chow et al. \cite{Chow}. The vector field $V$ generates the Ricci soliton viewed as a special solution of the Ricci flow. A Ricci soliton is said to be a gradient Ricci soliton, if $V = -\nabla f$ (up to a Killing vector field) for a smooth function $f$. Ricci solitons are also of interest to physicists who refer to them as quasi-Einstein metrics (for example, see Friedan \cite{Friedan}). \\

An odd dimensional analogue of Kaehler geometry is the Sasakian geometry. The Kaehler cone over a Sasakian Einstein manifold is a Calabi-Yau manifold which has application in physics in superstring theory based on a 10-dimensional manifold that is the product of the 4-dimensional space-time and a 6-dimensional Ricci-flat Kaehler (Calabi-Yau) manifold (see Candelas et al. \cite{Candelas}). Sasakian geometry has been extensively studied since its recently perceived relevance in string theory. Sasakian Einstein metrics have received a lot of attention in physics, for example, $p$-brane solutions in superstring theory, Maldacena conjecture (AdS/CFS duality) \cite{Maldacena}. For details, see Boyer, Galicki and Matzeu \cite{B-G-M}.\\

In \cite{Sharma} Sharma showed that if a $K$-contact (in particular, Sasakian) metric is a gradient Ricci soliton, then it is Einstein. This was also shown later independently by He and Zhu \cite{He-Zhu} for the Sasakian case. Recently, Sharma and Ghosh \cite{S-G} proved that a 3-dimensional Sasakian metric which is a non-trivial (i.e. non-Einstein) Ricci soliton, is homothetic to the standard Sasakian metric on $nil^3$. In this paper, we generalize these results and also answer the following question of H.-D. Cao (cited in \cite{He-Zhu}):``\textit{Does there exist a shrinking Ricci soliton on a Sasakian manifold, which is not Einstein?}", by proving
\begin{theorem}
If the metric of a ($2n+1$)-dimensional Sasakian manifold $M$ ($\eta, \xi, g, \varphi$) is a non-trivial (non-Einstein) Ricci soliton, then (i) $M$ is null $\eta$-Einstein (i.e. $D$-homothetically fixed and transverse Calabi-Yau), (ii) the Ricci soliton is expanding, and (iii) the generating vector field $V$ leaves the structure tensor $\varphi$ invariant, and is an infinitesimal contact $D$-homothetic transformation.
\end{theorem}
 
Conversely, we consider the following question: ``What can we say about an $\eta$-Einstein contact metric manifold $M$ which admits a vector field $V$ that leaves $\varphi$ invariant?" and answer it by assuming the invariance of the scalar curvature under $V$, in the form of the following result.

\begin{theorem}
If an $\eta$-Einstein contact metric manifold $M$ admits a vector field $V$ that leaves the structure tensor $\varphi$ and the scalar curvature invariant, then either $V$ is an infinitesimal automorphism, or $M$ is $D$-homothetically fixed and $K$-contact.
\end{theorem}

\begin{remark}Note that a Ricci soliton as a Sasakian metric is different from the Sasaki-Ricci soliton in the context of transverse Kaehler structure in a Sasakian manifold, for example see Futaki et al. \cite{F-O-W}).
\end{remark}

\begin{remark}Boyer et al. \cite{B-G-M} have studied $\eta$-Einstein geometry as a class of distinguished Riemannian metrics on contact metric manifolds, and proved the existence of $\eta$-Einstein metrics on many different compact manifolds. We would also like to point out that Zhang \cite{Zhang} showed that compact Sasakian manifolds with constant scalar curvature and satisfying certain positive curvature condition is $\eta$-Einstein.
\end{remark}

\begin{remark}Theorem 2 provides a generalization of the infinitesimal version of the following result of Tanno \cite{Tanno1} ``The group of all diffeomorphisms $\Phi$ which leave the structure tensor $\varphi$ of a contact metric manifold $M$ invariant, is a Lie transformation group, and coincides with the automorphism group $\mathcal{A}$ if $M$ is Einstein." Note that the scalar curvature of an Einstein metric is constant. We also note that the set of all vector fields on a contact metric manifold $M$, that leave $\varphi$ and scalar curvature invariant, forms a Lie sub-algebra of the Lie algebra of all smooth vector fields on $M$.
\end{remark}

\section{A Brief Review Of Contact Geometry}

A ($2n+1$)-dimensional smooth manifold is said to be contact if it has a global 1-form $\eta$ such that $\eta \wedge (d\eta)^{n} \neq 0$ on $M$. For a contact 1-form $\eta$ there exists a unique vector field $\xi$ such that $d\eta(\xi, X$)$ = 0$ and $\eta(\xi) = 1$. Polarizing $d\eta$ on the contact subbundle $\eta = 0$, we obtain a Riemannian metric $g$ and a (1,1)-tensor field $\varphi$ such that 
\begin{equation}\label{2}
d\eta(X,Y)=g(X,\varphi Y), \eta(X)=g(X,\xi), \varphi^{2}=-I+\eta \otimes \xi
\end{equation} 
$g$ is called an associated metric of $\eta$ and ($\varphi,\eta,\xi,g$) a contact metric structure. Following \cite{Blair} we recall two self-adjoint operators $h = \frac{1}{2}\pounds_\xi \varphi$ and $l = R(.,\xi)\xi$. The tensors $h$, $h\varphi $ are trace-free and $h\varphi = -\varphi h$ . We also have these formulas for a contact metric manifold.
\begin{equation}\label{3}
\nabla_{X}\xi = - \varphi X - \varphi hX
\end{equation}
\begin{equation}\label{4} 
\l -\varphi \l \varphi = -2(h^{2} + \varphi^{2})
\end{equation}\label{5}
\begin{equation}
\nabla_{\xi}h=\varphi-\varphi \l- \varphi h^{2}
\end{equation}
\begin{equation}\label{6}
Tr\l = Ric(\xi,\xi) = 2n - Tr h^{2}
\end{equation}
where $\nabla$, $R$, $Ric$ and $Q$ denote respectively, the Riemannian connection, curvature tensor, Ricci tensor and Ricci operator of $g$. For details see \cite{Blair}\\

A vector field $V$ on a contact metric manifold $M$ is said to be an infinitesimal contact transformation if $\pounds_V \eta = \sigma \eta$ for some smooth function $\sigma$ on $M$ . $V$ is said to be an infinitesimal automorphism of the contact metric structure if it leaves all the structure tensors $\eta, \xi, g, \varphi$ invariant (see Tanno \cite{Tanno2}).\\

A contact metric structure is said to be $K$-contact if $\xi$ is Killing with respect to $g$, equivalently, $h = 0$. The contact metric structure on $M$ is said to be Sasakian if the almost Kaehler structure on the cone manifold ($M \times R^{+}, r^2 g +dr^2$) over $M$, is Kaehler. Sasakian manifolds are $K$-contact and $K$-contact 3-manifolds are Sasakian. For a Sasakian manifold,
\begin{equation}\label{7}
(\nabla_{X}\varphi )Y = g(X, Y)\xi - \eta(Y)X
\end{equation}
\begin{equation}\label{8}
R(X,Y)\xi=\eta(Y)X-\eta(X)Y, \hskip.1in Q\xi = 2n \xi
\end{equation}
For a Sasakian manifold, the restriction of $\varphi$ to the contact sub-bundle $D$ ($\eta = 0$) is denoted by $J$ and ($D,J,d\eta$) defines a Kaehler metric on $D$, with the transverse Kaehler metric $g^T$ related to the Sasakian metric $g$ as $g = g^T + \eta \otimes \eta$. One finds by a direct computation that the transverse Ricci tensor $Ric^T$ of $g^T$ is given by
\begin{equation*}
Ric^T (X,Y)=Ric(X,Y)+2g(X,Y)
\end{equation*}
for arbitrary vector fields $X,Y$ in $D$. The Ricci form $\rho$ and transverse Ricci form $\rho^{T}$ are defined by
\begin{equation*}
\rho (X,Y)=Ric(X,\varphi Y), \hskip.1in \rho^{T}(X,Y)=Ric^{T} (X,\varphi Y)
\end{equation*}
for $X,Y \in D$. The basic first Chern class $2\pi c_{1}^B$ of $D$ is represented by $\rho^T$. In case $c_{1}^B = 0$, the Sasakian structure is said to be null (transverse Calabi-Yau). We refer to \cite{B-G-M} for details.\\

A contact metric manifold $M$ is said to be $\eta$-Einstein in the wider sense, if the Ricci tensor can be written as
\begin{equation}\label{9}
Ric(X,Y) = \alpha g(X,Y) + \beta \eta(X)\eta(Y)
\end{equation}
for some smooth functions $\alpha $ and $\beta $ on $M$. It is well-known (Yano and Kon \cite{Yano-Kon}) that $\alpha$ and $\beta$ are constant if $M$ is $K$-contact, and has dimension greater than 3.\\

Given a contact metric structure ($\eta, \xi, g, \varphi$), let $\bar{\eta} = a\eta, \bar{\xi}=\frac{1}{a}\xi, \bar{\varphi}=\varphi, \bar{g}=ag+a(a-1)\eta \otimes \eta$ for a positive constant $a$. Then ($\bar{\eta},\bar{\xi},\bar{\varphi},\bar{g}$) is again a contact metric structure. Such a change of structure is called a $D$-homothetic deformation, and preserves many basic properties like being $K$-contact (in particular, Sasakian). It is straightforward to verify that, under a $D$-homothetic deformation, a $K$-contact $\eta$-Einstein manifold transforms to a $K$-contact $\eta$-Einstein manifold such that $\bar{\alpha} =\frac{\alpha +2-2a}{a}$ and  $\bar{\beta}=2n-\bar{\alpha}$. We remark here that the particular value: $\alpha = -2$ remains fixed under a $D$-homothetic deformation, and as $\alpha + \beta = 2n$, $\beta$ also remains fixed. Thus, we state the following definition.
\begin{definition}
A $K$-contact $\eta$-Einstein manifold with $\alpha = -2$ is said to be $D$-homothetically fixed.
\end{definition}

\section{Proofs Of The Results}

\textbf{Proof Of Theorem 1: }Using the Ricci soliton equation (\ref{1}) in the commutation formula (Yano \cite{Yano}, p.23)
\begin{eqnarray}\label{10}
&&(\pounds_{V}\nabla_{X}g -\nabla_{X}\pounds_{V}g - \nabla_{[V,X]}g)(Y,Z)=\nonumber\\
&-&g((\pounds_{V}\nabla)(X, Y), Z) -g((\pounds_{V}\nabla)(X, Z), Y)
\end{eqnarray}
we derive
\begin{eqnarray}\label{11}
g((\pounds_{V}\nabla )(X, Y), Z) &=& (\nabla_{Z}Ric)(X, Y)\nonumber\\
&-&(\nabla_{X}Ric)(Y, Z) - (\nabla_{Y}Ric)(X, Z)
\end{eqnarray}
As $ \xi $ is Killing, we have $ \pounds_{\xi} Ric = 0 $ which, in view of (\ref{3}), the last equation of (\ref{8}) and $h = 0$, is equivalent to $ \nabla_{\xi} Q = Q\varphi - \varphi Q $. But for a Sasakian manifold, $ Q $ commutes with $ \varphi $, and hence $ Ric $ is parallel along $ \xi $. Moreover, differentiating the last equation of (\ref{8}), we have $(\nabla_{X}Q)\xi = Q\varphi X - 2n \varphi X $. Substituting $\xi $ for $Y$ in (\ref{11}) and using these consequences we obtain
\begin{equation}\label{12}
(\pounds_{V}\nabla)(X, \xi ) = - 2Q\varphi X + 4n \varphi X
\end{equation}
Differentiating this along an arbitrary vector field $Y $, using (\ref{7}) and the last equation of (\ref{8}), we find
\begin{eqnarray*}
(\nabla_{Y}\pounds_{V}\nabla )(X, \xi ) - (\pounds_{V}\nabla)(X, \varphi Y) = - 2(\nabla_{Y}Q)\varphi X + 2\eta(X)QY - 4n \eta(X)Y
\end{eqnarray*}
The use of the foregoing equation in the commutation formula \cite{Yano}:
\begin{equation}\label{13}
(\pounds_{V}R)(X, Y)Z = (\nabla_{X}\pounds_{V}\nabla)(Y, Z) - (\nabla_{Y}\pounds_{V}\nabla)(X, Z)
\end{equation} 
for a Riemannian manifold, shows that
\begin{eqnarray*}
(\pounds_{V}R)(X, Y)\xi - (\pounds_{V}\nabla)(Y, \varphi X) + (\pounds_{V}\nabla)(X, \varphi Y)= - 2(\nabla_{X}Q)\varphi Y \nonumber\\
+ 2(\nabla_{Y}Q)\varphi X + 2\eta(Y)QX - 2\eta(X)QY + 4n \eta(X)Y - 4n \eta(Y)X
\end{eqnarray*}
Substituting $ \xi $ for $ Y $ in the foregoing equation, using (\ref{12}) and the formula $ \nabla_{\xi}Q = 0 $ noted earlier, we find that 
\begin{equation}\label{14}
(\pounds_{V}R)(X, \xi)\xi = 4(QX - 2n X)
\end{equation}
Equation (\ref{1}) gives $ (\pounds_{V}g)(X, \xi) + 2(2n + \lambda )\eta(X) = 0$, which in turn, gives
\begin{equation}\label{15}
(\pounds_{V}\eta )(X) - g(\pounds_{V}\xi , X) + 2(\lambda + 2n)\eta(X) = 0 
\end{equation}
\begin{equation}\label{16}
\eta(\pounds_{V}\xi) = (2n + \lambda)
\end{equation} 
where we used the Lie-derivative of $g(\xi,\xi)=1$ along $V$. Next, Lie-differentiating the formula $ R(X, \xi)\xi = X - \eta(X)\xi $ [a consequence of the first formula in (\ref{8})] along $ V $, and using equations (\ref{14}) and (\ref{16}) provides
\begin{equation*}
4( QX - 2nX) - g(\pounds_{V}\xi, X)\xi + 2(2n + \lambda)X = -((\pounds_{V}\eta)(X))\xi 
\end{equation*}
By the direct application of (\ref{15}) to the the above equation we find
\begin{equation}\label{17}
 Ric(X, Y) = (n - \frac{\lambda}{2})g(X, Y) + (n + \frac{\lambda}{2})\eta(X)\eta(Y) 
\end{equation}
which shows that $M$ is $\eta$-Einstein with scalar curvature
\begin{equation}\label{18}
r = 2n(n + 1) - n\lambda
\end{equation}
At this point, we recall the following integrability formula \cite{Sharma}: 
\begin{equation}\label{19}
\pounds_{V}r = - \Delta r + 2\lambda r + 2  \left|Q\right|^{2}
\end{equation}
for a Ricci soliton, where $ \Delta r = -div Dr $. A straightforward computation using (\ref{17}) gives the squared norm of the Ricci operator as $ \left|Q\right|^{2} = 2n(n^{2} - n\lambda + \frac{\lambda^{2}}{4} + 4 n^{2}) $. Using this and (\ref{18}) in (\ref{19}), we obtain the quadratic equation $ (2n + \lambda )(2n + 4 - \lambda ) = 0 $. As $\lambda = -2n$  corresponds to $g$ becoming Einstein, we must have $ \lambda = 2n + 4 $ and hence the soliton is expanding, which proves part (ii). Moreover, equation (\ref{18}) reduces to $ r = -2n $. Thus equation (\ref{17}) assumes the form
\begin{equation}\label{20}
Ric(Y, Z) = -2 g(Y, Z) + 2(n + 1)\eta(Y)\eta(Z) 
\end{equation}
Hence, as defined in Section 2, $M$ is a $D$-homothetically fixed null $\eta$-Einstein manifold, proving part (i). Using (\ref{20}) in (\ref{11}) provides
\begin{equation}\label{21}
(\pounds_{V}\nabla)(Y, Z) = 4(n+1)\{\eta(Y)\varphi Z + \eta(Z)\varphi Y \}
\end{equation}
Differentiating this along $ X $, using equations (\ref{3}) and (\ref{7}), incorporating the resulting equation in (\ref{13}), and finally contracting at $X$ we get
\begin{equation}\label{22}
(\pounds_{V}Ric)(Y, Z) = 8(n+1)\{g(Y, Z) - (2n+1)\eta(Y)\eta(Z) \}
\end{equation}
Equation (\ref{20}) reduces the soliton equation (\ref{1}) to the form 
\begin{equation}\label{23}
(\pounds_{V}g)(Y, Z) = - 4(n+1)\{g(Y, Z) + \eta(Y)\eta(Z) \}
\end{equation}
Next, Lie-differentiating (\ref{20}) along $V$, and using (\ref{23}) shows
\begin{eqnarray}\label{24}
	(\pounds_{V}Ric)(Y, Z) &=& 8(n+1)\{g(Y, Z)+ \eta(Y)\eta(Z) \} \nonumber\\ 
	&+& 2(n+1)\{\eta(Z)(\pounds_{V}\eta)(Y) + \eta(Y)(\pounds_{V}\eta)Z \}
\end{eqnarray}
Comparing equations (\ref{22}) with (\ref{24}) and substituting $\xi $ for $Z$ leads to
\begin{equation}\label{25}
\pounds_{V}\eta=- 4(n+1)\eta
\end{equation}
Therefore, substituting $ \xi $ for $ Z $ in (\ref{23}) and using (\ref{25}) we immediately get $ \pounds_{V}\xi = 4 (n + 1)\xi $. Operating (\ref{25}) by $ d $, noting $ d$ commutes with $\pounds_{V}$ and using the first equation of (\ref{2}) we find 
\begin{equation*}
(\pounds_{V}d\eta )(X, Y) = - 4( n + 1)g(X, \varphi Y)
\end{equation*}
Its comparison with the Lie-derivative of the first equation of (\ref{2}) and the use of (\ref{23}) yields $ \pounds_{V}\varphi = 0 $, completing the proof.\\

Before proving Theorem 2, we state and prove the following lemma.
\begin{lemma}
If a vector field $V$ leaves the structure tensor $\varphi$ of the contact metric manifold $M$ invariant, then there exists a constant $ c $ such that\\
(i)$\pounds_{V}\eta = c\eta $, (ii)$ \pounds_{V}\xi = - c \xi $, (iii) $ \pounds_{V}g = c(g + \eta \otimes \eta ) $.
\end{lemma}

\noindent
Though this lemma was proved by Mizusawa in \cite{Mizusawa}, to make the paper self-contained, we provide a slightly different proof as follows.\\

\noindent
\textbf{Proof: }Lie-differentiating the formulas $\varphi \xi = 0$ and $\eta (\varphi X) = 0$ and using $\pounds_V \varphi = 0$, we find $\pounds_V \xi = -c \xi$, and $\pounds_V \eta = c \eta$ for a smooth function $c$ on $M$. Next, Lie-derivative of the formula $\eta(X)=g(X,\xi)$ along $V$ gives
\begin{equation}\label{26}
(\pounds_V g)(X,\xi)=2c\eta(X)
\end{equation}
The Lie-derivative of the first equation of (\ref{2}) along $V$ provides
\begin{equation}\label{27}
(\pounds_V g)(X,\varphi Y)=((dc)\wedge \eta)(X,Y)+cg(X,\varphi Y)
\end{equation}
Substituting $\xi$ for $Y$ in the above equation we get $dc = (\xi c)\eta$. Taking its exterior derivative, and then exterior product with $\eta$ shows $(\xi c)(d\eta)\wedge \eta = 0$. By definition of the contact structure, $(d\eta)\wedge \eta$ is nowhere zero on $M$, and so $\xi c = 0$. Hence $dc = 0$, i.e. $c$ is constant. Using this consequence, and equations (\ref{26}) and (\ref{27}) we obtain (iii), completing the proof.\\

\noindent
\textbf{Proof Of Theorem 2 : }By virtue of Lemma 1, we have 
\begin{equation}\label{28}
(\pounds_{V}g)(Y, Z) = c\{g(Y, Z) + \eta (Y) \eta (Z) \}
\end{equation}
Differentiating this and using (\ref{3}) we get
\begin{equation}\label{29}
(\nabla_{X}\pounds_{V}g)(Y,Z)= - c\{\eta(Z)g(Y,\varphi X + \varphi hX) + \eta(Y)g(Z,\varphi X + \varphi hX) \}
\end{equation}
Equation (\ref{10}) can be written
\begin{equation}\label{30}
(\nabla_{X}\pounds_{V}g)(Y,Z)= g((\pounds_{V}\nabla)(X,Y),Z) + g((\pounds_{V}\nabla)(X,Z),Y)
\end{equation}
A straightforward computation using (\ref{29}) and (\ref{30}) shows
\begin{equation*}
(\pounds_{V}\nabla)(Y,Z) = - c \{\eta(Z)\varphi Y + \eta(Y)\varphi Z + g(Y, \varphi hZ)\xi\}
\end{equation*}
Its covariant differentiation and  use of (\ref{5}) provides
\begin{eqnarray*}
(\nabla_{X}\pounds_{V}\nabla)(Y,Z) &=& - c \{\eta(Z)(\nabla_{X}\varphi)Y + \eta(Y)(\nabla_{X}\varphi)Z \nonumber\\
&-& g(Z,\varphi X + \varphi hX)\varphi Y - g(Y,\varphi X + \varphi hX)\varphi Z  \nonumber\\
&-& g(\varphi hY, Z)(\varphi X + \varphi hX) + g((\nabla_{X}\varphi h)Y, Z)\xi\}
\end{eqnarray*}  
Using this in the commutation formula (\ref{13}) for a Riemannian manifold, contracting at $X$, and using equations (\ref{2}), (\ref{3}) and also the well known formula: $(div \varphi )X = -2n\eta(X)$ for a contact metric (see \cite{Blair}), we find 
\begin{eqnarray}\label{31}
(\pounds_{V}Ric)(Y, Z) &=& c \{- 2 g(Y,Z) + 2 g(hY, Z) \nonumber\\ &+& 2(2n+1)\eta(Y)\eta(Z) \}
 - c g((\nabla_{\xi}\varphi h)Y, Z)
\end{eqnarray}
Also, Lie-differentiating (\ref{9}) along $V$ and using Lemma 1 we have 
\begin{equation}\label{32}
(\pounds_V Ric)(Y,Z) = (V\alpha + c\alpha) g(Y,Z) + (V\beta + c(\alpha + 2\beta)) \eta(Y)\eta(Z)
\end{equation}
Comparing the previous two equations shows that 
\begin{eqnarray*}
	[V\alpha + c(\alpha + 2)]g(Y, Z) + [V\beta +c(\{\alpha + 2\beta - 2(2n + 1)\}] \eta(Y)\eta(Z) \nonumber\\
	-c[ 2 g(hY, Z) - g((\nabla_{\xi}\varphi h)Y, Z)] = 0
\end{eqnarray*}
On one hand, we substitute $Y=Z=\xi$ in the above equation getting one equation, and on the other hand, we contract the above equation (noting that both $h$ and $\varphi h$ are trace-free) getting another equation. Solving the two equations we obtain
\begin{equation}\label{33}
V\alpha + c(\alpha +2)=0, \hskip.1in V\beta +c(\alpha +2\beta-4n-2)=0
\end{equation}
The $g$-trace of equation (\ref{9} ) gives the scalar curvature
\begin{equation}\label{34}
r=(2n+1)\alpha + \beta
\end{equation}
The divergence of (\ref{9}) along with the contracted second Bianchi identity yields $dr = 2d\alpha +2(\xi \beta)\eta$. Taking its exterior derivative, and then exterior product with $\eta$ we have $(\xi \beta)\eta \wedge d\eta = 0$. As $\eta \wedge d\eta$ vanishes nowhere on $M$, we find $\xi \beta = 0$ whence $dr = 2d\alpha$. Hence $V\alpha = Vr = 0$, by hypothesis. Thus, it follows from (\ref{34}) that $V\beta = 0$. Consequently, equations (\ref{33}) reduce to: $c(\alpha + 2) = 0$ and $c(\alpha + 2\beta -4n-2)=0$, and hence imply that, either $c = 0$ in which case $V$ is an infinitesimal automorphism, or $\alpha = -2$ and $\alpha + 2\beta = 4n+2$. In the second case, adding the two equations gives $\alpha + \beta = 2n$. But, from equation (\ref{9}) we have $\alpha + \beta = Tr. l$. Therefore, $Tr. l = 2n$, and applying equation (\ref{6}) we obtain $h = 0$, i.e. $M$ is $K$-contact. As $\alpha = -2$, the $\eta$-Einstein structure is $D$-homothetically fixed, completing the proof.

\section{An Explicit Example}An explicit example of non-trivial Ricci soliton as a Sasakian metric is the (2n+1)-dimensional Heisenberg group $\mathcal{H}^{2n+1}$ (which arose from quantum mechanics) of matrices of type $
\left[
\begin{array}{ccc}
1	&	Y		&	z\\
O^t	&	I_n		&	X^t\\
0	&	O		&	1
\end{array}
\right]
$, where $X=(x_1,...,x_n), Y=(y_1,...,y_n), O=(0,...,0) \in R^n, z \in R$. As a manifold, this is just $R^{2n+1}$ with coordinates ($x^{i},y^{i},z$) where $i = 1,...,n$, and has the left-invariant Sasakian structure ($\eta, \xi, \varphi, g$) defined by $\eta = \frac{1}{2}(dz- \sum_{i=1}^n y^i dx^i)$, $\xi = 2\frac{\partial}{\partial z}$, $\varphi (\frac{\partial}{\partial x^i})=-\frac{\partial}{\partial y^i}$, $\varphi (\frac{\partial}{\partial y^i})=\frac{\partial}{\partial x^i}+y^i \frac{\partial}{\partial z}$, $\varphi (\frac{\partial}{\partial z})=0$, and the Riemannian metric $g = \eta \otimes \eta + \frac{1}{4}\sum_{i=1}^{n}((dx^i)^2 + (dy^i)^2)$. Its $\varphi$-sectional curvature (i.e. the sectional curvature of plane sections orthogonal to $\xi$) is equal to $-3$, so its Ricci tensor satisfies equation (\ref{20}), as shown by Okumura \cite{Okumura}, and hence $\mathcal{H}^{2n+1}$ is a $D$-homothetically fixed null $\eta$-Einstein manifold. Setting $V = \sum_{i=1}^{n}(V^{i}\frac{\partial}{\partial x^i}+ \bar{V}^{i}\frac{\partial}{\partial y^i})+V^{z}\frac{\partial}{\partial z}$, using equations: $\pounds_V \xi = 4(n+1)\xi$, $\pounds_V \varphi=0$ obtained in the proof of Theorem 1, and the aforementioned actions of $\varphi$ on the coordinate basis vectors, shows that $V^i$ and $\bar{V}^i$ do not depend on $z$ and yields the PDEs:
\begin{eqnarray*}
\frac{\partial V^i}{\partial x^j}&=&\frac{\partial \bar{V}^i}{\partial y^j},\hskip .1in \frac{\partial V^i}{\partial y^j}=-\frac{\partial \bar{V}^i}{\partial x^j}, \hskip .1in y^{i}\frac{\partial V^i}{\partial y^j}=\frac{\partial V^z}{\partial y^j}\nonumber\\
\bar{V}^j &=&y^{j}\frac{\partial V^z}{\partial z}-y^{i}\frac{\partial \bar{V}^i}{\partial y^j}, \hskip .1in \frac{\partial V^z}{\partial z}=-4(n+1)
\end{eqnarray*}
The last equation readily integrates as $V^z =-4(n+1)z+F(x^i ,y^i)$. For a special solution, assuming $F=0$, $V^i =cx^i$, $\bar{V}^i =cy^i$ and substituting in the above PDEs, we get $c=-2(n+1)$, and hence the Ricci soliton vector field $V=-2(n+1)(x^i \frac{\partial}{\partial x^i}+y^i \frac{\partial}{\partial y^i}+2z\frac{\partial}{\partial z})$. For dimension 3, this reduces to $V=-4(x\frac{\partial}{\partial x}+y\frac{\partial}{\partial y}+2z\frac{\partial}{\partial z})$ which occurs on p. 37 of \cite{Chow} without the factor 4, but gets adjusted with our $\lambda = 6$ which is 4 times their $\lambda = 3/2$.

\begin{remark} Another conclusion that we draw for Theorem 1 is the following: The value $-2n$ for the scalar curvature $r$ obtained during the proof, and the equation (\ref{17}) show that the generalized Tanaka-Webster scalar curvature \cite{Blair} $ W = r - Ric(\xi, \xi) + 4n$ vanishes.
\end{remark}
\textbf{Acknowledgment:} We thank the referee for valuable suggestions. R.S. was supported by University Research Scholar grant. This work is dedicated to Bhagawan Sri Sathya Sai Baba and Sri Ramakrishna Paramahansa.

\end{document}